\documentclass{svmult}   % please, do NOT add options!
\usepackage{graphicx}    % standard LaTeX graphics tool when
\usepackage{amsmath,amsfonts,amssymb,epsfig,epic}
\usepackage{graphicx}
\newtheorem{thm}{Theorem}
\newtheorem{defn}{Definition}

\newcommand{\be}[1]{\begin{equation}\label{#1}}
\newcommand{\ee}{\end{equation}}

\begin{document}

%%%%%%%%%%%%%%%%%   TITLE OF THE PAPER   %%%%%%%%%%%%%%%%%%%%%%%

\title*{Almost Global Convergence in Singular Perturbations of Strongly
  Monotone Systems} % Put your title here

\titlerunning{Singular Perturbations of Strongly Monotone Systems}

%%%%%%%%%%%%%%%%%       AUTHOR NAMES     %%%%%%%%%%%%%%%%%%%%%%%

\author{
Liming Wang\inst{1}
\and Eduardo D. Sontag\inst{2}}

%%%%%%%%%%%%%%%%%   AUTHOR NAMES FOR INDEX GENERATION   %%%%%%%%%%%%

\index{Wang L.}
\index{Sontag E.D.}

%%%%%%%%%%%%%%%%%         AFFILIATIONS       %%%%%%%%%%%%%%%%%%%%%%%

\institute{Department of Mathematics, Rutgers University, New Brunswick, NJ
  08903, USA, 
\texttt{wshwlm@math.rutgers.edu}
\and 
Department of Mathematics, Rutgers University, New Brunswick, NJ 08903, USA,
\texttt{sontag@math.rutgers.edu}
}

\maketitle
\renewcommand{\abstractname}{Abstract.}

\begin{abstract}

This paper deals with global convergence to equilibria, and in particular
Hirsch's generic convergence theorem for strongly monotone systems, for
singular perturbations of monotone systems.

\end{abstract}

\section{Introduction}
Monotone systems constitute a rich class of models for which global and
almost-global convergence properties can be established.  They are
%eds3: adding some CDC references to show ``relevance''
particularly useful in biochemical models (see discussion and references
in~\cite{04sysbio,05ejc}), and also appear in areas like coordination
(\cite{moreau}) and other problems in control theory (\cite{chisci}).
This paper studies extensions, using geometric singular perturbation theory,
of Hirsch's generic convergence theorem for monotone systems
(\cite{Hirsch2,Hirsch,Hirsch-Smith,Smith}).
%%%%%%%%%%%%%%%%%%%%%%%%%%%%%%%%%%%%%%%%%%%%%%%%%%%%%%%%%%%%%%%%%%%%%%%%%
Informally stated, Hirsch's result says that almost every bounded solution of
a strongly monotone system converges to the set of equilibria.  There is a rich
literature regarding the application of this powerful theorem, as well as of
other results dealing with everywhere convergence when equilibria are unique
(\cite{Smith,Dancer,JiFa}), to models of biochemical systems.
Unfortunately, many models in biology are not monotone.
%%%%%%%%%%%%%%%%%%%%%%%%%%%%%%%%%%%%%%%%%%%%%%%%%%%%%%%%%%%%%%%%%%%%%%%%%
In order to address this drawback (as well as to study properties of large
systems which are monotone but which are hard to analyze in their entirety),
a recent line of work introduced an input/output approach that is based on the
analysis of interconnections of monotone systems.
%%%%%%%%%%%%%%%%%%%%%%%%%%%%%%%%%%%%%%%%%%%%%%%%%%%%%%%%%%%%%%%%%%%%%%%%%
For example, the approach allows one to view a  \emph{non-}monotone system as a
``negative'' feedback loop of monotone open-loop systems, thus leading to
results on global stability and the emergence of oscillations under
transmission delays, and to the construction of relaxation oscillators by slow
adaptation rules on feedback gains.  See~\cite{04sysbio,05ejc} for expositions
and many references.
% for longer paper: add from below list!
The present paper is in the same character.

Our motivation arose from the observation that time-scale separation may also
lead to monotonicity.  This point of view is of special interest in the 
context of biochemical systems; for example, Michaelis Menten kinetics are
mathematically justified as singularly perturbed versions of mass action
kinetics.
%refer to Keshet, Murray, 
A system that is not monotone may become monotone once that fast variables are
replaced by their steady-state values.
%%%%%%%%%%%%%%%%%%%%%%%%%%%%%%%%%%%%%%%%%%%%%%%%%%%%%%%%%%%%%%%%%%%%%%%%%
A trivial linear example that illustrates this point is
$\dot x$$\,=\,$$-x$$-y$, $\varepsilon \dot y$$\,=\,$$-y$$+$$x$, with $\varepsilon $$>$$0$.
This system is not monotone with respect to any orthant cone.
On the other hand, for $\varepsilon \ll1$, the fast variable $y$ tracks $x$,
so the slow dynamics is well-approximated by $\dot x=-2x$
(which is strongly monotone, because every scalar system is).

We consider systems $\dot x=f(x,y)$, $\varepsilon \dot y=g(x,y)$ for which the reduced system
$\dot x=f(x,h(x))$ is strongly monotone (in fact, a slightly stronger technical condition on derivatives is assumed) and the fast system $\dot y=g(x,y)$ has a
%eds3: comma was separated from $x$:
unique globally asymptotically stable steady state $y=h(x)$ for each $x$, and satisfies an input to state stability type of property with
respect to $x$.
One may expect that the original system inherits global (generic) convergence
properties, at least for all $\varepsilon $$>$$0$ small enough, and this is indeed 
the object of our study.  
%%%%%%%%%%%%%%%%%%%%%%%%%%%%%%%%%%%%%%%%%%%%%%%%%%%%%%%%%%%%%%%%%%%%%%%%%
This question may be approached in several ways.  One may view $y-h(x)$
as an input to the slow system, and appeal to the theory of asymptotically
autonomous systems.  Another approach, the one that we develop here, is
through geometric invariant manifold theory 
(\cite{Fenichel,Jones,Nipp}).
There is a manifold $M_\varepsilon $, invariant for the full dynamics, which attracts
all near-enough solutions, with an asymptotic phase property.
The system restricted to the invariant manifold $M_\varepsilon $ is a regular
perturbation of the fast ($\varepsilon $$=$$0$) system.  As remarked in Theorem 1.2 in
Hirsch's early paper~\cite{Hirsch2}, a $C^1$ regular perturbation
of a flow with eventually positive derivatives also has generic convergence.
So, solutions in the manifold will be generally well-behaved, and asymptotic
phase implies that solutions track solutions in $M_\varepsilon $, and hence also
converge to equilibria if solutions on $M_\varepsilon $ do.
A key technical detail is to establish that the tracking solutions also start
from the ``good'' set of initial conditions, for generic solutions of the
large system.  

For simplicity, we discuss here only the case of cooperative systems
(monotonicity with respect to the main orthant), but proofs in the
case of general cones are similar and will be discussed in a paper
under preparation.

\addtolength{\textheight}{-3cm}   % This command serves to balance the column lengths

%%%%%%%%%%%%%%%%%%%%%%%%%%%%%%%%%%%%%%%%%%%%%%%%%%%%%%%%%%%%%%%%%%%%%%%%%%%%%%%%
\section{Statement of Main Result}

We are interested in systems in singularly perturbed form:
\begin{eqnarray}
\frac{dx}{dt}&=&f(x,y) \label{eqn:general} \\
\varepsilon \frac{dy}{dt}&=&g(x,y),\nonumber
\end{eqnarray}
where $x \in \mathbb R^n$, $y \in \mathbb R^m$, $0<\varepsilon \ll 1$, 
and $f$ and $g$ are smooth functions.
We will present some preliminary results in general, but
for our main theorem we will restrict attention to the case when $g$ has the
special form $g(x,y)=Ay+h(x)$, where $A$ is a Hurwitz matrix (all eigenvalues
have negative real part) and $h$ is a smooth function.
That is, we will specialize to systems of the following form:
\begin{eqnarray}
\frac{dx}{dt}&=&f(x,y) \label{eqn:full} \\
\varepsilon \frac{dy}{dt}&=&Ay+h(x). \nonumber
\end{eqnarray}
(We remark later how our results may be extended to a broader class of systems.)
%In the fast time scale $\tau=t/\varepsilon$, (\ref{eqn:full}) is governed by:
%\begin{eqnarray}
%\frac{dx}{d\tau}&=&\varepsilon f(x,y) \label{eqn:fast} \\
%\frac{dy}{d\tau}&=&Ay+h(x). \nonumber
%\end{eqnarray}
%(A more general situation is when this form is only attained locally, after a diffeomorphism.)
%In the slow time scale, setting $\varepsilon$ to zero, we have:
Setting $\varepsilon$ to zero, we have:
\begin{eqnarray}
\frac{dx}{dt}&=&f(x,m_0(x)), \label{eqn:limit} 
\end{eqnarray}
where $m_0(x)=-A^{-1}h(x)$. 
As usual in singular perturbation theory, our goal is to use properties of the
limiting system (\ref{eqn:limit}) in order to derive conclusions about the
full system (\ref{eqn:full}) when $0<\varepsilon \ll 1$. 
%lm8:
In this paper, $A \subset B$ means that $A$ is a strict subset of $B$, while $A \subseteq B$ contains the case of $A=B$.
%lm4 We will assume given two sets $K$ and $L$ 
%eds3: comma in US-style English when listing items, before ``and''
% see: http://www.getitwriteonline.com/archive/021201.htm
We will assume given three sets $K$, $\widetilde K$, and $L$
which satisfy the following hypotheses
%lm7
(some technical terms are defined later):
\begin{description}
\item{\bf H1} The set 
%lm4:$K$
$\widetilde K$ is an $n$-dimensional $C^\infty$ simply connected compact manifold with boundary.

\item{\bf H2} The set $L$ is a bounded open subset of $\mathbb R^m$,
and 
%lm9:
\[
M_0=\{(x,y)\ |\ y=m_0(x),\  x \in \widetilde K\},
\]
the graph of $m_0$, is contained in $\widetilde K \times L$.

\item{\bf H3} The flow $\{\psi_t\}$ of the limiting system (\ref{eqn:limit}) has eventually positive derivatives on $\widetilde K$.

\item{\bf H4} The set 
%lm7
$\widetilde K$ is convex, and therefore it is p-convex too.

\item{\bf H5} For each $\varepsilon>0$ sufficiently small, the forward trajectory
  under (\ref{eqn:full}) of each point in 
%lm4: $D=\mbox{Int}K \times L$ is precompact in $D$.
$\widetilde D=\mbox{Int} \widetilde K \times L$ is precompact in $\widetilde D$.

\item{\bf H6} The equilibrium set 
%lm4  $E_0=\{x \in \mbox{Int}K \ | \  f(x,m_0(x))=0\}$ is countable. 
$E_0=\{x \in \mbox{Int}\widetilde K \ | \  f(x,m_0(x))=0\}$ is countable. 

%lm4: add another assumption, the original H6 becomes H7
%lm5: change to Int\widetilde K
\item{\bf H7} The set $K \subset \mbox{Int}\widetilde K$ is compact, and for each $\varepsilon>0$ sufficiently small, the set $D=K \times L$ is positively invariant.

\end{description}

Note that the equilibria of (\ref{eqn:full}) do not depend on $\varepsilon$, and 
%lm4 they
the ones in $\widetilde D$ 
are in 1-1 correspondence with elements of $E_0$.
The main theorem is:
\begin{thm}
Under assumptions {\bf H1-H7}, 
%lm4 for each sufficiently small $\varepsilon>0$, the
there exists $\varepsilon^*>0$ such that for each $0<\varepsilon<\varepsilon^*$, the
forward trajectory of (\ref{eqn:full}) starting from almost every point in $D$
converges to some equilibrium. 
\label{thm:main}
\end{thm}

\noindent{\bf Remark:}
A variant of this result is to assume that the reduced system (\ref{eqn:limit})
has a unique equilibrium.  In this case, one may improve the conclusions of
the 
%lm5 Theorem
theorem to global (not just generic) convergence, by appealing to
results of Hirsch and others that apply when equilibria are unique.
The proof is simpler in that case, since the foliation structure given by
Fenichel's theory (see below) is not required.
%see also, for future work:
In the opposite direction, one could drop the assumption of countability and
instead provide theorems on generic convergence to the set of equilibria, or
even to equilibria if hyperbolicity conditions are satisfied, in the 
spirit of what is done in the theory of strongly monotone systems.
% strong mon needed so that epsilon-systems are monotone, it seems.. so not
% ji-fa or dancer...
\qed

\section{Terminology}

The following standard terminology is defined for a general ordinary
differential equation: 
\begin{eqnarray}
\frac{d z}{dt}=F(z), \label{eqn:def}
\end{eqnarray}
where $F:\mathbb R^N \rightarrow \mathbb R^N$ is a $C^1$ vector field. For any $z \in \mathbb R^{N}$, we denote the maximally defined solution of (\ref{eqn:def}) with initial condition $z$ by $t \rightarrow \phi_t(z), t \in I(z)$, where $I(z)$ is an open interval in $\mathbb R$ that contains zero. For each $t \in \mathbb R$, the set of $z \in \mathbb R^N$ for which $\phi_t(z)$ is defined is an open set $W(t) \subseteq \mathbb R^N$, and $\phi_t: W(t) \rightarrow W(-t)$ is a diffeomorphism.
The collection of maps $\{\phi_t \}_{t \in R}$ is called the flow of (\ref{eqn:def}).
We also write just $z(t)$ for the solution of (\ref{eqn:def}), if the initial
condition $z(0)$ is clear from the context. 
The forward trajectory of $z \in \mathbb R^{N}$ is a parametrized curve $t \rightarrow \phi_t(z)$. Its image is the forward orbit of $z$, denoted as $O_+(z)$. The backward trajectory and the backward orbit $O_-(z)$ are defined analogously.
A set $U \subseteq \mathbb R^{N}$ is \emph{positively (respectively, negatively) invariant} if $O_+(U) \subseteq U$. It is \emph{invariant} if it is both positively and negatively invariant.

We borrow the notation from \cite{Jones} for the forward evolution of a set $U \subseteq V \subseteq \mathbb R^N$ restricted to $V$:
\[
U\cdot_V t=\{\phi_t(p): \ p \in U \mbox{ and } \phi_s(p) \in V \mbox{ for all } 0 \leq s \leq t\}.
\]
Let us denote the interior and the closure of a set $U$ as Int$U$ and $\overline U$ respectively.

\begin{defn}
The flow $\{\phi_t\}$ of (\ref{eqn:def}) is said to have \emph{eventually
  positive derivatives on a set $V \subseteq \mathbb R^N$} if there exists
  $t_0$ such that $\frac{\partial \phi_t^i}{\partial z_j}(z)>0$ for all $t
  \geq t_0$, $z \in V$. 
\label{def:EPD}
\end{defn}
%lm8:add remark
%lm9:grammar
If (\ref{eqn:def}) is of dimension one, i.e. $N=1$, $\{\phi_t\}$ has eventually positive derivatives automatically. In practice, the following sufficient condition is easier to check. If the vector field of (\ref{eqn:def}) satisfies $\frac{\partial F_i}{\partial z_j}(z) \geq 0$, for all $z \in V$, $i \not= j$, and the matrix $\frac{\partial F}{\partial z}(z)$ is irreducible for all $z \in V$, then $\{\phi_t\}$ has eventually positive derivatives. (This condition is not necessary.)

\begin{defn}
An open set $W \subseteq \mathbb R^N$ is called \emph{p-convex} if $W$
contains the entire line segment joining $x$ and $y$ whenever $x,y \in V$ and
$x \leq y$, where $x \leq y$ means $x_i \leq y_i$ for all $i=1, \cdots, N$. 
\end{defn}
%lm4: add a lemma, Hirsch thm 4.4
The next lemma is a restatement of theorem 4.4 in \cite{Hirsch}:
\begin{lemma}
%eds3:
Suppose that
the open set $W \subseteq R^n$ is  $p$-convex and the flow $\{\phi_t\}$ of (\ref{eqn:def}) has eventually positive derivatives on $W$. Let $W^c\subseteq W$ be the set of points whose forward orbit has compact closure in $W$. If the set of equilibrium points is countable, then $z(t)$ converges to a equilibrium as $t \rightarrow \infty$, for almost every $z \in W^c$. 
\label{lemma:Hirsch_Con}
\end{lemma}
%lm4
The following fact follows from differentiability of solutions with respect to
``regular'' perturbations in the dynamics; see \cite{Hirsch}, Theorem 1.2:

\begin{lemma}
Assume $V \subset W$ is a compact set in which the flow $\{\phi_t\}$ has eventually positive derivatives. Then, there exists $\delta>0$ with the following property. Let $\{\psi_t\}$ denote the flow of a $C^1$ vector field $G$ such that the $C^1$ norm of $F(z)-G(z)$ is less than $\delta$ for all $z$ in $V$. Then there exists $t_*>0$ such that if $t\geq t_*$ and $\psi_s(z) \in V$ for all $s \in [0,t]$, then $\frac{\partial \psi_t^i}{\partial z^j}(z)>0$.
\label{lemma:Hirsch}
\end{lemma}

The Appendix reviews
% (definitions \ref{def:half_space} and \ref{def:manifd_bdy}) 
the definition of a ``$C^r$ ($1 \leq r \leq\infty$) manifold $M$ with
boundary'' in the sense used on geometric singular perturbation theory
%%%%%%%%%%%%%%%%%%
We denote the boundary of such a manifold $M$ as
$\partial M$, and denote $M \setminus \partial M$ as Int$M$, when there is no
confusion with the notation for the interior of a set. 

\begin{defn}
A compact, connected $C^r$ manifold $M \subset \mathbb R^N$ with boundary is
said to be \emph{locally invariant}
 under the flow of (\ref{eqn:general}) if for each
$p$ in Int$M$ there exists a time interval 
$I_p=(t_1,t_2)$, for some $t_1<0<t_2$,
such that $\phi_t(p) \in M$ for all $t \in I_p$.  
\end{defn}
%lm5 added the fast system, since the def of NH make sense only in the fast system.
When $\varepsilon \not= 0$, we can ``stretch'' the time ($\tau=\frac{t}{\varepsilon}$), and consider the fast system:
\begin{eqnarray}
\frac{dx}{d\tau}&=&\varepsilon f(x,y) \label{eqn:f_general} \\
\frac{dy}{d\tau}&=&g(x,y).\nonumber
\end{eqnarray}
The above system is equivalent to (\ref{eqn:general}). The corresponding fast system for (\ref{eqn:full}) is
\begin{eqnarray}
\frac{dx}{d\tau}&=&\varepsilon f(x,y) \label{eqn:f_full} \\
\frac{dy}{d\tau}&=&Ay+h(x).\nonumber
\end{eqnarray}

\begin{defn}
%lm3:
Let $M$ be an $n$-dimensional manifold (possibly with boundary) contained in $\{(x,y) \ |\ g(x,y)=0\}$. We say that $M$ is normally hyperbolic relative to (\ref{eqn:f_general}) if all eigenvalues of the matrix $\frac{\partial g}{\partial y}(p)$ have nonzero real part for every $p \in M$.
\end{defn}

\section {Proof of the Main Theorem}

Recall the definition of 
%lm4: $M_0=\{(x,y)\,|\,\ y=m_0(x), x \in K\}$. Since $K$ is
$M_0=\{(x,y)\,|\,\ y=m_0(x), x \in \widetilde K\}$. Since $\widetilde K$ is
an $n$-dimensional $C^\infty$ compact manifold with boundary, and $m_0$ is a
smooth function, $M_0$ is also an $n$-dimensional $C^\infty$ compact manifold
with boundary. 

Our proofs are based on Fenichel's theorems~\cite{Fenichel}, in the forms
presented and developed by Jones in~\cite{Jones}.

\begin{figure}
  \centering \includegraphics[scale=0.7,angle=0]{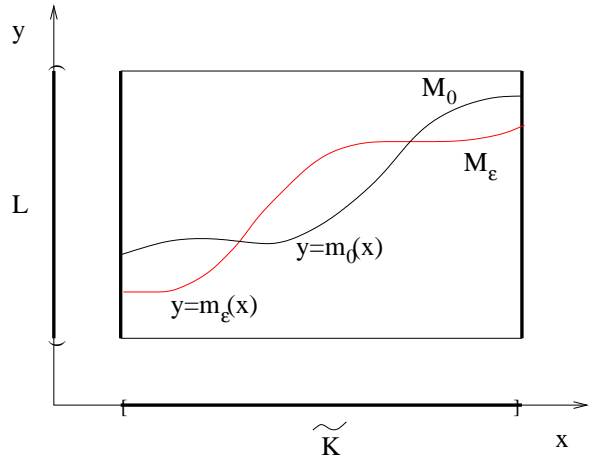}
  \caption{For simplicity, we sketch manifolds $M_\varepsilon$ and $M_0$ of a system where $n=m=1$. 
The set $\widetilde K$ is a compact set in $x$, and $L$ is an open set in $y$. The red curve denotes the locally invariant manifold $M_\varepsilon$ and the black curve denotes $M_0$.}\label{fig:manifold}
\end{figure}

\medskip
\noindent{\bf Fenichel's First Theorem}
{\em Under assumption H1, if $M_0$ is normally hyperbolic relative to (\ref{eqn:f_full}), then there exists $\varepsilon_0>0$, such that for every $0< \varepsilon<\varepsilon_0$ and $r>0$, 
there is a function $y=m_\varepsilon(x)$, defined on 
$\widetilde K$, of class $C^r$ jointly in $x$ and $\varepsilon$, such that 
\[
M_\varepsilon=\{(x,y)\ |\ y=m_\varepsilon(x),\  x \in \widetilde K\}
\]
is locally invariant under (\ref{eqn:full}), see Figure \ref{fig:manifold}.}

The requirement that $M_0$ be normally hyperbolic is satisfied in our case, as $g(x,y)=Ay+h(x)$ and therefore $\frac{\partial g}{\partial y}(p)=A$, which is 
%lm16:
Hurwitz, for each $p \in M_0$. 

We will pick a particular $r>1$ in the above theorem from now on. 

Let us interpret local invariance in terms of equations. Let $(x(t),y(t))$ be the solution to (\ref{eqn:full}) with initial condition $(x_0,y_0)$, such that  
%lm4: $x_0 \in \mbox{Int}K$ 
$x_0 \in \mbox{Int}\widetilde K$ and $y_0=m_\varepsilon(x_0)$. Local invariance implies that $(x(t),y(t))$ satisfies
\begin{eqnarray}
\frac{dx(t)}{dt}&=& f\big(x(t),m_\varepsilon(x(t))\big) \label{eqn:epsilon} \\
y(t)&=&m_\varepsilon(x(t)), \label{eqn:yx}
\end{eqnarray}
for all $t$ small enough. Actually, this is also true for all $t \geq 0$. The
argument is as follows. By 
%lm4 {\bf H3}
{\bf H5}, $(x(t),y(t))$ is 
well-defined and remains in 
%lm4 $D$
$\widetilde D$ for all
$t \geq 0$. Let $T=\{t \geq 0 \ | \ y(t)=m_\varepsilon(x(t))\}$. Then, $T$ is not
empty, and $T$ is closed by the continuity of $m_\varepsilon(x(t))$ and
$y(t)$. Also, $T$ is open, since $M_\varepsilon$ is locally invariant. So $T=\{t
\geq 0\}$, that is, $x(t)$ is a solution to (\ref{eqn:epsilon}) and
$y(t)=m_\varepsilon(x(t))$ for all $t \geq 0$. 

In (\ref{eqn:epsilon}), the $x$-equation is decoupled from the $y$-equation, which allows us to reduce to studying a lower-dimension system. Another advantage is that, as $\varepsilon$ approaches zero, the limit of system (\ref{eqn:epsilon}) is system (\ref{eqn:limit}), which describes the flows on $M_0$. If $M_0$ has some desirable property, it is natural to expect this property is inherited by the perturbed manifold $M_\varepsilon$. An example of this principle is provided by the following lemma.
\begin{lemma}
%lm8: H1-H4
Under assumptions {\bf H1}-{\bf H3}, for each $0<\varepsilon<\varepsilon_0$, the flow $\{\psi_t\}$ of (\ref{eqn:epsilon}) has eventually positive derivatives on Int$\widetilde K$.
\label{lemma:positive}
\end{lemma}

Proof. Applying Lemma~\ref{lemma:Hirsch}, there exist $\delta>0$ such that when the $C^1$ norm of $m_0(x)-m_\varepsilon(x)$ is less than $\delta$ for all $x \in \widetilde K$, there exists $t_*>0$ with the property that: if $t\geq t_*$ and $\psi_s(x) \in \widetilde K$ for all $s \in [0,t]$, then $\frac{\partial \psi^i_t}{\partial x^j}(x)>0$. Since $m_\varepsilon$ is of class $C^r$, jointly in $x$ and $\varepsilon$, we can pick $\varepsilon>0$ small enough to control $\delta$. If we then prove Int$\widetilde K$ is invariant under (\ref{eqn:epsilon}), that is, for any $x_0 \in \mbox{Int}\widetilde K$, the solution $x(t)$ of (\ref{eqn:epsilon}) with initial condition $x_0$ stays in Int$\widetilde K$ for all $t \geq 0$, then we are done. Let us now prove that Int$\widetilde K$ is positively invariant under (\ref{eqn:epsilon}). 

Let $y_0=m_\varepsilon(x_0)$ and $y(t)=m_\varepsilon(x(t))$. Then, $(x(t),y(t))$ is the solution to (\ref{eqn:full}) with initial condition $(x_0,y_0) \in \widetilde D$. By {\bf H5},  $(x(t),y(t))$ stays in $\widetilde D$ for all $t \geq 0$, and therefore $x(t) \in \mbox{Int}\widetilde K$ for all $t \geq 0$.
\qed

Flows with eventually positive derivatives have particularly appealing properties, as in Lemma \ref{lemma:Hirsch_Con}. To apply that lemma, we need to check two conditions. First, for every point in Int$\widetilde K$, its forward trajectory under (\ref{eqn:epsilon}) has compact closure in Int$\widetilde K$. Second, the number of equilibria of (\ref{eqn:epsilon}) is countable. Suppose that the first property does not hold, and let $x(t)$ be a solution to (\ref{eqn:epsilon}) with $x(0) \in \mbox{Int}\widetilde K$ but $\lim_{j \rightarrow \infty} x(t_j) \notin \mbox{Int}\widetilde K$ for some sequence $\{t_j\}$. So, $\big(x(t),m_\varepsilon(x(t))\big)$ is a solution for (\ref{eqn:full}), and its forward orbit is not precompact in $\widetilde D$. This violates {\bf H5}. To check the second condition, we introduce the following sets: 
\[
E_\varepsilon=\{x \in \mbox{Int}\widetilde K \ | \ f(x,m_\varepsilon(x))=0\}\ ,\ 0 \leq \varepsilon \ll 1.
\]
We claim that $E_\varepsilon \subseteq E_0$ for all $\varepsilon$ small enough, which implies that $E_\varepsilon$ is countable, by {\bf H6}.
 Let us prove the claim. It is clear that $E_0$ consists of the $x$-coordinates of all equilibria of (\ref{eqn:full}) in $\widetilde D$. Fix $0<\varepsilon<\varepsilon_0$ and pick any $x_0 \in E_\varepsilon$, $y_0=m_\varepsilon(x_0)$. The solution $(x(t),y(t))$ to (\ref{eqn:full}) with initial condition $(x_0,y_0)$ satisfies (\ref{eqn:epsilon}) and (\ref{eqn:yx}) for $t$ small enough. But 
\[
\frac{dx(t)}{dt}=f \big( x(t),y(t))=f \big(x(t), m_\varepsilon(x(t))\big) \equiv 0,
\]
so $(x(t),y(t))=(x_0,y_0)$ for all $t \geq 0$, and therefore $x_0 \in E_0$.
Applying Lemma \ref{lemma:Hirsch_Con} we have:

\begin{lemma}
Under assumptions {\bf H1}-{\bf H6}, for each $0<\varepsilon<\varepsilon_0$, there exists a set $\mathcal C_\varepsilon \subseteq \mbox{Int}\widetilde K$ such that the forward trajectory of (\ref{eqn:epsilon}) for every point of $\mathcal C_\varepsilon$ converges to some equilibrium, and the measure of Int$\widetilde K \setminus \mathcal C_\varepsilon$ is zero.
\label{lemma:converge}
\end{lemma}

Until now, we have discussed the flow only when restricted to the locally
invariant manifold $M_\varepsilon$. The next theorem, stated in the form given by~\cite{Jones}, deals with more global behavior. In \cite{Jones}, the theorem is stated for $\varepsilon>0$, but some properties also hold for $\varepsilon=0$ (\cite{Jones_pc}). (We will apply this result again with a fixed $r>1$.) The notation $[-\delta,\delta]$ stands for the cube $\{\ (y_1, \dots, y_m) \ | \ -\delta \leq |y_i| \leq \delta\ \}$.

\medskip
\noindent{\bf Fenichel's Third Theorem}
{\em Let $\varepsilon_0$ be as in Fenichel's First Theorem. Under assumption H1, if $M_0$ is normally hyperbolic relative to (\ref{eqn:f_full}), then there exists $0 < \varepsilon_1<\varepsilon_0$
%lm8: remove $\delta_0$ and all $\delta$s from functions h and W^s etc.
and $\delta>0$ such that for every $0 \leq \varepsilon<\varepsilon_1$ and $r>0$, there is a function  
\[
%lm7
h_\varepsilon:\widetilde K \times[-\delta,\delta]\rightarrow \mathbb R^n
\]
%(\bar p, \lambda)$, defined for all $\bar p \in
%M_\varepsilon$ and $|\lambda| \leq \delta$, of class $C^r$ jointly in $\bar p$
%and $\varepsilon$, 
such that the following properties hold: 
\begin{enumerate}
\item 
%lm4 For each $\bar p\in M_\varepsilon$, $h_\varepsilon(\bar p,0)=\bar x$, where $\bar x$ is the $x$-coordinate of $\bar p$. 
%lm7
For each $x \in \widetilde K$, $h_\varepsilon(x,0)=x$.
\item 
The image of the map 
\begin{eqnarray*}
T_\varepsilon: \widetilde K \times [-\delta,\delta] &\rightarrow& \mathbb R^n \times \mathbb R^m \\
(x,\lambda) &\mapsto& ( h_\varepsilon(x,\lambda), \lambda + m_\varepsilon( h_\varepsilon(x,\lambda) ))
\end{eqnarray*}
is defined as the stable manifold $W^s_\varepsilon(M_\varepsilon)$ of $M_\varepsilon$. 
%lm8:
For $p=(x, m_\varepsilon(x)) \in M_\varepsilon$, the stable fibers $W^s_\varepsilon(p)$, defined as $T_\varepsilon(\{x\} \times [-\delta,\delta])$, form a ``positively invariant'' family when $\varepsilon \not=0$, in the sense that
\[
W^s_\varepsilon(p) 
%eds2: made scripts bigger so can see better
\cdot_{W^{\scriptstyle s}_{\scriptstyle\varepsilon}(M_{\scriptstyle\varepsilon})} t 
\subseteq W^s_\varepsilon(\phi_t(p)).
\]
%lm7
\item ``Asymptotic Phase''. There are positive constants $k$ and $\alpha$ such that for any $p, q \in \mathbb R^{n+m}$, if $q \in W^s_\varepsilon(p)$, $\varepsilon \not=0$, then
\[
|\phi_t(p)-\phi_t(q)| \leq ke^{-\alpha t} 
\]
for all $t \geq 0$ as long as $\phi_t(p)$ and $\phi_t(q)$ stay in $W^s_\varepsilon(M_\varepsilon)$.

\item The stable fibers are disjoint, i.e., for $q_i \in W^s_\varepsilon(p_i)$, $i=1,2$, either $W^s_\varepsilon(p_1)\bigcap W^s_\varepsilon(p_2)=\emptyset$ or $W^s_\varepsilon(p_1)=W^s_\varepsilon(p_2)$.

\item The function $h_\varepsilon(x, \lambda)$ is $C^r$ jointly in $\varepsilon$, $x$ and $\lambda$. When $\varepsilon=0$, $h_{0,\delta}(x, \lambda)=~x$.
\end{enumerate}
}
%lm8
\begin{figure}
  \centering \includegraphics[scale=0.55,angle=0]{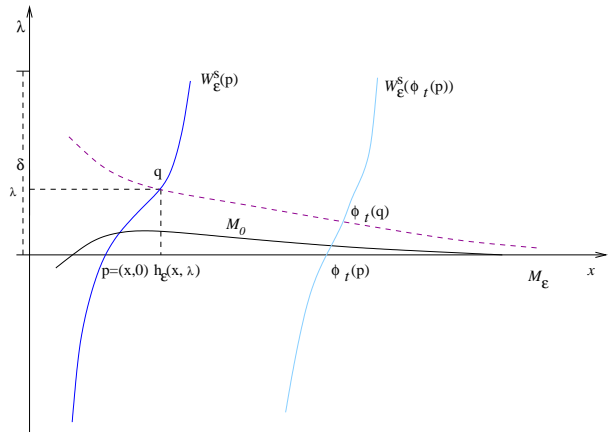}
\caption{To illustrate the geometric meaning of Fenichel's Third Theorem, we
  sketch the locally invariant manifold and stable fibers of a system, in the
  case of $n$$=$$m$$=$$1$. The dimensions of the manifolds $M_\varepsilon$, $M_0$, and  stable fibers are one. $M_\varepsilon$ is the graph of $\lambda=0$, and $M_0$ is the graph of $m_0(x)-m_\varepsilon(x)$ (black
  curve). These manifolds may intersect at some equilibrium points.
  Through each point $p \in M_\varepsilon$ ($x$-axis),
  there is a stable fiber $W^s_\varepsilon(p)$ (blue curve).
  We call $p$ the ``base point'' of the
  fiber. The fiber consists of the pairs
  $(h_\varepsilon(x,\lambda),\lambda)$, where $|\lambda|\leq \delta$.
  If a solution (purple dashed curve)
  starts on fiber $W^s_\varepsilon(p)$, after a small time $t$, it
  evolves to a point on another stable fiber
  $W^s_\varepsilon(\phi_t(p))$ (light blue curve); 
   this is the ``positive invariance'' property.} 
\label{fig:fiber}
\end{figure}

The next lemma gives a sufficient condition to guarantee that a point is on some fiber.
\begin{lemma}
%lm7
Let $\varepsilon_1$ and $\delta$ be as in Fenichel's Third Theorem. There exists $0<\varepsilon_2<\varepsilon_1$, such that for every $0<\varepsilon<\varepsilon_2$, the set 
\[
\mathcal A_\delta:=\{(x,y) \ | \ x \in K, \ |y-m_0(x)| \leq \frac{\delta}{2}\}
\]
is a subset of $W^s_\varepsilon(M_\varepsilon)$.
\label{lemma:band}
\end{lemma}
To prove this lemma, we need the following result:
\begin{lemma}
Let $U$ and $V$ be compact, convex subsets of $\mathbb R^n$ and $\mathbb R^m$
respectively. Suppose given a continous function
\begin{eqnarray*}
\phi:\ U \times V &\rightarrow& \mathbb R^n \times \mathbb R^m \\
(x,y) &\mapsto& (\phi_1(x,y),\phi_2(x,y))
\end{eqnarray*}
satisfying $\|\phi_1(x,y)-x\| \leq \rho_1$, $\|\phi_2(x,y)-y\| \leq \rho_2$ for some $\rho_1>0$, $\rho_2>0$ and all $(x,y) \in U \times V$. Then every point $(\alpha, \beta) \in U \times V$ with dist$(\alpha,\partial U) \geq \rho_1$ and dist$(\beta, \partial V) \geq \rho_2$ is in the image of $\phi$.
\label{lemma:close_id}
\end{lemma}
Proof. For such a point $(\alpha, \beta)$, consider the map $\Phi(x,y)=(\Phi_1(x,y),\Phi_2(x,y)):=(x,y)-(\phi_1(x,y), \phi_2(x,y))+(\alpha,\beta)$. Thus $\Phi$ maps $U \times V$ into itself. If not, say $\Phi_1(x,y)$ is not in $U$, that is, $x-\phi_1(x,y)+\alpha$ is not in $U$. Since $\|x-\phi_1(x,y)\| \leq \rho_1$, then dist$(\alpha, \partial U)<\rho_1$, contradiction. The case when $\Phi_2(x,y)$ is not in $V$ follows similarly. Since $\Phi$ maps $U \times V$ into itself, and the product of convex sets is still convex, by Brouwer's Fixed Point Theorem, there is some $(\bar x, \bar y)\in U \times V$ so that $\Phi(\bar x, \bar y)=(\bar x, \bar y)$, which means that $(\phi_1(\bar x, \bar y), \phi_2(\bar x, \bar y))=(\alpha, \beta)$, as we wanted to prove. 
\qed

%lm8 changed as we discussed last time.
Proof of Lemma \ref{lemma:band}. Define for each $0 \leq \varepsilon<\varepsilon_1$, the map
\begin{eqnarray*}
\phi_\varepsilon: \ \widetilde K \times [-\delta, \delta] &\rightarrow& \mathbb R^n \times \mathbb R^m \\
(x,\lambda) &\mapsto& (h_\varepsilon(x,\lambda), \lambda+m_\varepsilon(h_\varepsilon(x,\lambda))-m_0(h_\varepsilon(x,\lambda))).
\end{eqnarray*}

By property 5 in Fenichel's Third Theorem and the compactness of $\widetilde K$, we have 
\[
\|h_\varepsilon((x,m_\varepsilon(x)),\lambda)-x\| \leq C_1(\varepsilon),
\]
and
\[
\|[\lambda+m_\varepsilon(h_\varepsilon(x,\lambda))-m_0(h_\varepsilon(x,\lambda))]-\lambda\| \leq C_2(\varepsilon)
\]
for some positive functions $C_i$ of $\varepsilon$ such that $C_i \rightarrow 0$ as $\varepsilon \rightarrow 0$, $i=1,2$. For any such $\varepsilon$, apply Lemma \ref{lemma:close_id} with $U=\widetilde K$, $V=[-\delta,\delta]$, $\rho_1=C_1(\varepsilon)$, $\rho_2=C_2(\varepsilon)$ and $\phi=\phi_\varepsilon$. Since dist$(K, \partial \widetilde K)$ and $\frac{\delta}{2}$ are independent of $\varepsilon$, we can pick $0<\varepsilon_2<\varepsilon_1$ such that dist$(K, \partial \widetilde K) > C_1(\varepsilon)$ and $\frac{\delta}{2}>C_2(\varepsilon)$ for all $\varepsilon \in (0,\varepsilon_2)$. By Lemma~\ref{lemma:close_id}, 
\begin{equation}
K \times [-\frac{\delta}{2},\frac{\delta}{2}] \subset \phi_{\varepsilon}(\widetilde K \times [-\delta,\delta])
\label{eqn:subset}
\end{equation}

Define a map
\begin{eqnarray*}
\pi_0: \ \mathbb R^n \times \mathbb R^m &\rightarrow& \mathbb R^n \times \mathbb R^m \\
(x,\lambda) &\mapsto& (x,\lambda+m_0(x))
\end{eqnarray*}

Consider the composition of $\pi_0$ and $\phi_\varepsilon$: By property 2 in Fenichel's Third Theorem, its image, $\pi_0\circ \phi_\varepsilon(\widetilde K \times [-\delta, \delta])$, is the stable manifold $W^s_\varepsilon(M_\varepsilon)$ of $M_\varepsilon$. According to (\ref{eqn:subset}), $\pi_0(K \times [-\frac{\delta}{2}, \frac{\delta}{2}]) \subset W^s_\varepsilon(M_\varepsilon)$. Notice that $\mathcal A_\delta=\pi_0(K \times [-\frac{\delta}{2}, \frac{\delta}{2}])$, we are done.
\qed

\begin{lemma}
Let $\varepsilon_2$ be as in Lemma \ref{lemma:band}, and $\delta$ as in Fenichel's Third Theorem. Under assumption {\bf H7}, there exists $0<\varepsilon_3<\varepsilon_2$ such that for each $0<\varepsilon<\varepsilon_3$, if $p \in D$, then there exists $T_0 \geq 0$, and $\phi_t(p) \in \mathcal A_{\delta}$ for all $t\geq T_0$. 
\label{lemma:attract}
\end{lemma}
Proof. Setting $z=y-m_0(x)$ and $\tau=\frac{t}{\varepsilon}$, (\ref{eqn:full}) becomes
\begin{eqnarray*}
\frac{dx}{d\tau}&=&\varepsilon f(x,z+m_0(x)) \\
\frac{dz}{d\tau}&=&Az-\varepsilon m_0'(x)f(x,z+m_0(x)).
\end{eqnarray*}
So 
\[
z(\tau)=z(0)e^{A\tau}-\varepsilon \int_0^\tau e^{A(\tau-s)} m_0'(x)f\big( x,z+m_0(x)\big) \, ds.
\]
Notice that $\|e^{At}\|\leq Ce^{\beta t}$, for some positive constant $C$ and negative $\beta$, which is greater than the real part of all eigenvalues of $A$. So,
\[
\left \| \varepsilon \int_0^\tau e^{A(\tau-s)} m_0'(x)f\big( x,z+m_0(x)\big) \, ds \right \| \leq \frac{2\varepsilon M C}{|\beta|},
\]
where $M$ is an upper bound of the function $\|m_0'(x)f(x,y)\|$ on $\overline D$.
Let
\[
\varepsilon=\frac{\delta |\beta|}{8MC} \quad\mbox{and}\quad
T_0'=\max \{\frac{1}{|\beta|}\ln\frac{4C\|z(0)\|}{\delta},0\} \,.
\]
Then, we have
$\|z(\tau)\| \leq \frac{\delta}{2}$
for all $\tau \geq T_0'$. Back to the slow time scale, we let $T_0=\varepsilon T_0'$. Therefore, $\phi_t(p) \in \mathcal A_{\delta}$ for all $t \geq T_0$,
%eds3: remind me please why we need to know that we stay in IntK?
%lm6: If x is not in K, then the point may not be on any fiber. The point of having those added lemmas is to prove there is a band with x in K not \widetilde K.if the $x$-coordinate of $\phi_t(p)$ stays in $K$, but this can be easily
derived from {\bf H7}.

\medskip
\noindent{\bf Remark:}
Except for the normal hyperbolicity assumption, 
Lemma \ref{lemma:attract} is the only place where the special
structure~(\ref{eqn:full}) was used.
Consider a more general system as in (\ref{eqn:general}), and assume that
$g(x,m_0(x))=0$ on $\widetilde K$ for some smooth function $m_0$. By the same change of
variables as in the above proof, (\ref{eqn:general}) is equivalent to  
\begin{eqnarray*}
\frac{dx}{d\tau}&=&\varepsilon f(x,z+m_0(x)) \\
\frac{dz}{d\tau}&=&g(x,z+m_0(x))-\varepsilon m_0'(x)f(x,z+m_0(x)).
\end{eqnarray*}
The only property that we need in the lemma is that for any initial condition
$(x(0),z(0))$, the solution $(x(t),z(t))$ satisfies 
\[
\limsup_{t\rightarrow\infty} |z(t)| \leq \gamma\left(\limsup_{t\rightarrow\infty} d(t)\right)
\]
where $\gamma$ is a function of class $\cal K$, that is to say, a continuous
function $[0,\infty)\rightarrow[0,\infty)$ with $\gamma(0)=0$, and
$d(t)=\varepsilon m_0'(x(t))f\big(x(t),z(t)+m_0(x(t))\big)$. 
In terms of the functions $m_0$ and $g$, we may introduce the control system
$dz/dt = G(d(t),z)+u(t)$, where $d$ is a compact-valued
``disturbance'' function and $u$ is an input, and $G(d,z)=g(d,z+m_0(d))$.
Then, the property of input-to-state stability with input $u$ (uniformly on
$d$), which can be characterized in several different manners, including by
means of Lyapunov functions,
%TODO to add refs for journal paper
%(see~\cite{TAC89,yuanTACnewchars,cime04})
provides the desired condition.
\qed

Lemma \ref{lemma:attract} proves that every trajectory in $D$ is attracted to 
$A_{\delta}$ and therefore is also attracted to $M_\varepsilon$. This will lead to our proof of the main theorem.

\subsubsection*{Proof of the main theorem.}
Choose $\varepsilon^*=\varepsilon_3$, defined in Lemma \ref{lemma:attract}. For any $p \in D$, there are three cases:
\begin{enumerate}
\item $p \in M_\varepsilon$. By Lemma \ref{lemma:converge}, the forward trajectory converges to an equilibrium except for a set of measure zero.
\item $p \in \mathcal A_{\delta} \subset W^s_\varepsilon(M_\varepsilon)$. 
Then $p$ is on some fiber, say $W^s_\varepsilon(\bar p)$, where $\bar p=(\bar x, m_\varepsilon(\bar x)) \in M_\varepsilon$. 
If $\bar x$ is in $\mathcal C_\varepsilon$ (defined in Lemma \ref{lemma:converge}), then
$\phi_t(\bar p) \rightarrow q$, for some $q \in E_0$. By the ``asymptotic phase'' property of Fenichel's Third Theorem, $\phi_t(p)$ also converges to $q$. To deal with the case when $\bar x \not\in \mathcal C_\varepsilon$, it is enough to show that the set
\[
\mathcal B_\varepsilon=\bigcup_{\bar x \in \mbox{Int}\widetilde K \setminus \mathcal C_\varepsilon} W^s_{\varepsilon}(\bar p)
\]
as a subset of $\mathbb R^{m+n}$ has measure zero. 
Define
\[
\mathcal F_\varepsilon = \left(\mbox{Int}\widetilde K\setminus
\mathcal C_\varepsilon\right)\times[-\delta,\delta].
\]
Since $\mbox{Int}\widetilde K\setminus \mathcal C_\varepsilon$
has measure zero in $\mathbb R^n$, also $\mathcal F_\varepsilon$ has
measure zero. On the other hand $T_\varepsilon(\mathcal F_\varepsilon)=\mathcal B_\varepsilon$, and Lipschitz maps send measure zero sets to measure zero sets, we are done.

\item $p \in D \setminus \mathcal A_{\delta}$. By Lemma
\ref{lemma:attract}, $\phi_t(p) \in \mathcal A_{\delta}$ for all $t
\geq T_0$. Without loss of generality, we assume that $T_0$ is an integer. If
$\phi_{T_0}(p) \in \mathcal A_{\delta} \setminus \mathcal
B_\varepsilon$, then $\phi_t(p)$ converges to an equilibrium. Otherwise,
$p \in \bigcup_{k \geq 0, k \in \mathbb Z} \phi_{-k} (\mathcal
B_\varepsilon)$. 
Since the set $\mathcal B_\varepsilon$ has measure zero and
$\phi_{-k}$ is Lipschitz, $\phi_{-k}(\mathcal B_\varepsilon)$ has measure zero for all $k$, and the countable union of them still has measure zero.
\qed
\end{enumerate}

\section{An Example}
%lm8: revised the example
%lm9: consider a more general example, y is now a vector
%lm11: generalized more by Prof. Sontag, x is now a vector
Consider the following system:
\begin{eqnarray}
\frac{dx_i}{dt}&=&\gamma_i(y_1,\dots,y_m)-\beta_i(x_1, \dots, x_n), \ \ i=1, \dots, n, \nonumber \\
\varepsilon\frac{dy_j}{dt}&=&-d_jy_j-\alpha_j(x_1, \dots, x_n), \ \ d_j>0, \ \ j=1, \dots, m,
\label{eqn:example}
\end{eqnarray}
where $\alpha_j$, $\beta_i$ and $\gamma_i$ are smooth functions. We assume that
%lm10:remove assumption that \beta is odd, and add assumptions of smoothness
%lm11:change \beta and \gamma to vector function 
\begin{enumerate}
\item 
%lm15: changed conditions
%When $n>1$, for all $i,k=1, \dots, n$ and $j,l=1,\dots, m$, 
%eds4: added:
%and all $x\in\mathbb R^n$,
%the partial derivatives
%\[
%\frac{\partial \beta_i}{\partial x_k}(x)<0 \mbox{ for } i \not=k, 
%\ \ 
%\frac{\partial \alpha_j}{\partial x_k}(x)\geq 0, 
%\ \ \frac{\partial \gamma_i}{\partial y_l}(x)\leq 0.
%\]
%lm11: I separate the case of n=1 and n>1, because in our example of limit cycle, d \gamma / dy > 0, but we do need \gamma or \alpha to have negative partial derivatives for n>1.
%eds4: yes, but we should eventually give a unified condition (not now!)
%      also, I am changing slightly this:
%While for $n=1$, there is no assumptions on the partial derivatives.
%(For the case $n=1$, no conditions are imposed on the partial derivatives.)
%lm11: I only \beta_i approaches \pm \infty when x_i does.
%eds5: this was not right, I think.  One needs to be careful about what
%      happens for the remaining variables even as x_i is large.
%\item As $x_i$ approaches $\pm \infty$, $\beta_i$ approaches $\pm \infty$. 
The reduced system
\[
\frac{d x_i}{dt}= \gamma_i(-\frac{\alpha_1}{d_1}, \dots, -\frac{\alpha_m}{d_m})-\beta_i(x_1, \dots, x_n):=F_i(x_1, \dots, x_n), \ \ i=1,\dots, n
\]
%eds8: put equation numbers in next three
has partial derivatives that satisfy:
\be{reducedpartialexample}
\frac{\partial F_i}{\partial x_k}=\sum_{l=1}^m -\frac{1}{d_i}\frac{\partial \gamma_i}{\partial y_l}\frac{\partial \alpha_l}{\partial x_k}-\frac{\partial \beta_i}{\partial x_k}>0 \mbox{ for } i \not=k.
\ee
\item
%lm15 changed back
%lm12:change conditions.
%The function $\beta_i$ satisfies
%\begin{eqnarray}
%\lim_{x_1 \rightarrow +\infty, \dots, x_n \rightarrow +\infty} \beta_i(x_1,\dots,x_n)=+\infty \nonumber \\
%\lim_{x_1 \rightarrow -\infty, \dots, x_n \rightarrow -\infty} \beta_i(x_1,\dots,x_n)=-\infty.
% \label{cond:infty}
%\end{eqnarray}
For each $i$,
\be{betacond1}
\lim_{u\rightarrow+\infty} \min_{x\in S_{\scriptstyle i}(u)} \beta_i(x) = +\infty
\ee
and
\be{betacond2}
\lim_{u\rightarrow-\infty} \max_{x\in S_{\scriptstyle i}(u)} \beta_i(x) = -\infty
\ee
where $S_i(u)$ is the set of vectors in $\mathbb R^n$ whose $i$th coordinate
is $u$.
(For $n=1$, this means simply that 
$\lim_{x\rightarrow\pm\infty} \beta_i(x) = \pm\infty$.)
\item There exists a positive contant $M_j$ such that $|\alpha_j(x)| \leq M_j$ for all $x \in \mathbb R^n$.
\item The number of roots of the 
%eds4
%set 
system
of equations
\[
\gamma_i(\alpha_1(x), \dots, \alpha_m(x))=\beta_i(x), \ \ i=1,\dots, m
\]
is countable.
\end{enumerate}
%lm15: remove
%lm12: 
%Since $\frac{\partial \beta_i}{\partial x_k}(x)<0$ for all $i \not=k$, (\ref{cond:infty}) implies 
%\[
%\lim_{u\rightarrow+\infty} \min_{x\in S_{\scriptstyle i}(u)} \beta_i(x) = +\infty
%\]
%\[
%\lim_{u\rightarrow-\infty} \max_{x\in S_{\scriptstyle i}(u)} \beta_i(x) = -\infty,
%\]
%where $S_i(u)$ is the set of vectors in $\mathbb R^n$ whose $i$th coordinate
%is $u$.
%eds5:
%The conditions are very natural.  The condition on the $\beta_i$'s is
%satisfied, for example, if there is a linear decay term $-x_i$ in the
%differential equation for $x_i$, and all other variables appear saturated in
%this rate.

We are going to show that on any large enough region, and provided that
$\varepsilon$ is sufficiently small, almost every trajectory converges to an equilibrium. To emphasize the need for small $\varepsilon$, we also show that
when $\varepsilon>1$, a limit cycle could appear.

%TODO: prove this, for the journal paper?:
%It is easy to verify that the origin is the only equilibrium, and the Jacobian of
%the right-hand side of~(\ref{eqn:example}) at the origin
%has trace $1-\frac{1}{\varepsilon}$ and determinant $\frac{1}{\varepsilon}$.  Thus,
%when $0<\varepsilon<1$, the trace is negative and the determinant is
%positive, so $(0,0)$ is locally stable.
%We wish to obtain global conclusions.

To apply our main theorem, we take
%lm9:
%\[
%K=[-a,a],\ \widetilde K=[-a-1, a+1],\ L=(-b,b),\ D=K \times L,\ \widetilde D=\mbox{Int}\widetilde K \times L,
%\]
\[
L=\{\,y \in \mathbb R^m \ | \ |y_j| <b_j, \ j =1, \dots, m \,\},
\]
where $b_j$ is an arbitrary positive number greater than $\frac{M_j}{d_j}$. Picking such $b_j$ assures $y_j \frac{d y_j}{d t}<0$ for all $x \in \mathbb R$ and $|y_j|=b_j$, i.e. the vector field points transverally inside on the boundary of $L$. Let
%lm10: change the definition of K, a_1, and a_2
%lm11: change again
\[
K=\{\, x \in \mathbb R^n \ | \ -a_{i,2} \leq x_i \leq a_{i,1}, \ i=1,\dots, n\,\}
\]
where $a_{i,1}$ and $a_{i,2}$ can be any positive numbers such that 
%eds5:
\[
%\beta(a_{i,1})>N_i:=\max_{|y_i| \leq b_i} |\gamma_i(y_1, \dots, y_m)|, \ \ \beta(-a_{i,2})< -N_i.
\beta_i(x)>N_i:=\max_{|y_j| \leq b_j} |\gamma_i(y_1, \dots, y_m)|
\]
whenever $x \in \mathbb R^n$ satisfies that its $i$th coordinate
$x_i\geq a_{i,1}$, and 
\[
\beta_i(x)< -N_i
\]
whenever $x \in \mathbb R^n$ satisfies that its $i$th coordinate
$x_i\leq -a_{i,2}$.
All large enough $a_{i,j}$'s satisfy this condition, because of the assumption
made on $\beta$.

So, we have $x_i\frac{d x_i}{d t}<0$ for all $y \in L$, $x_i=a_{i,1}$ and $x_i=-a_{i,2}$. We then take 
\[
\widetilde K=\{\, x \in \mathbb R^n \ | \ -a_{i,2}-1 \leq x_i \leq a_{i,1}+1,\ i=1,\dots n \,\},
\]
$D=K \times L$ and $\widetilde D=\mbox{Int}\widetilde K \times L$. Thus, the vector field will point into the interior of $D$ and $\widetilde D$. Hypotheses {\bf H5} and {\bf H7} follow directly from this fact. (Sketch: {\bf H7} is obvious. Suppose {\bf H5} does not hold.  Then, there exists some solution $(x(t), y(t))$ of (\ref{eqn:example}) in $D$, and a sequence $\{t_j\}$ such that $(x(t_j), y(t_j)) \rightarrow (\bar x, \bar y)$ as $j \rightarrow \infty$. Suppose that $\bar y_k=b_k$ for some $k \in \{\,1, \dots, m\,\}$, and $\{y_k({t_j})\}$ is strictly increasing to $b_k$. This will contradict the fact that
$\frac{dy_i}{dt}<0$ above the $y_k$-nullcline. The other cases follow
similarly.) 
%lm15: remove
%lm11:
%Now let us verify {\bf H3}. If $n=1$, the reduced system is of dimension one, so {\bf H3} holds. When $n>1$, the reduced system is 
%\[
%\frac{d x_i}{dt}= \gamma_i(-\frac{\alpha_1}{d_1}, \dots, -\frac{\alpha_m}{d_m})-\beta_i(x_1, \dots, x_n):=F_i(x_1, \dots, x_n), \ \ i=1,\dots, n.
%\]
%The partial derivatives
%\[
%\frac{\partial F_i}{\partial x_k}=\sum_{l=1}^m -\frac{1}{d_i}\frac{\partial \gamma_i}{\partial y_l}\frac{\partial \alpha_l}{\partial x_k}-\frac{\partial \beta_i}{\partial x_k}>0 \mbox{ for } i \not=k.
%\]
%So, the reduced system has eventually positive derivatives on $\widetilde K$ by our remark after Definition \ref{def:EPD}.
{\bf H3} follows from our assumption 1, and it is easy to see the other hypotheses also hold. By our main theorem, for sufficiently small $\varepsilon$, the forward trajectory of (\ref{eqn:example}) starting from almost every point in $D$ converges to some equilibrium.

%lm8: (In this example, the limiting system $\frac{dx}{dt}=(-4\frac{e^x-1}{e^x+1})^m-\frac{1}{3}x^3+x$ is strongly monotone
%and has a unique equilibrium.  In this case, one may prove that the forward
%trajectory starting from every, not merely a generic, point in $D$ converges
%to the origin.)

On the other hand, convergence does not hold for large $\varepsilon$. Let
\[
n=1, \ \beta_1(x_1)=\frac{x_1^3}{3}-x_1, \ m=1, \ \alpha_1(x_1)=4\tanh x_1,\  \gamma(y_1)=y_1, \ d_1=1.
\]
It is easy to verify that $(0,0)$ is the only equilibrium. When 
%lm10:
$\varepsilon>1$, the trace of the Jacobian at $(0,0)$ is
%lm13: determinant is 15/epsilon not 1/epsilon
$1-\frac{1}{\varepsilon}>0$, its determinant is $\frac{15}{\varepsilon}>0$,
so the (only) equilibrium in $D$ is repelling. By the Poincar\'e-Bendixson Theorem, there exists a limit cycle in $D$.

\begin{figure}
  \centering \includegraphics[scale=0.4,angle=0]{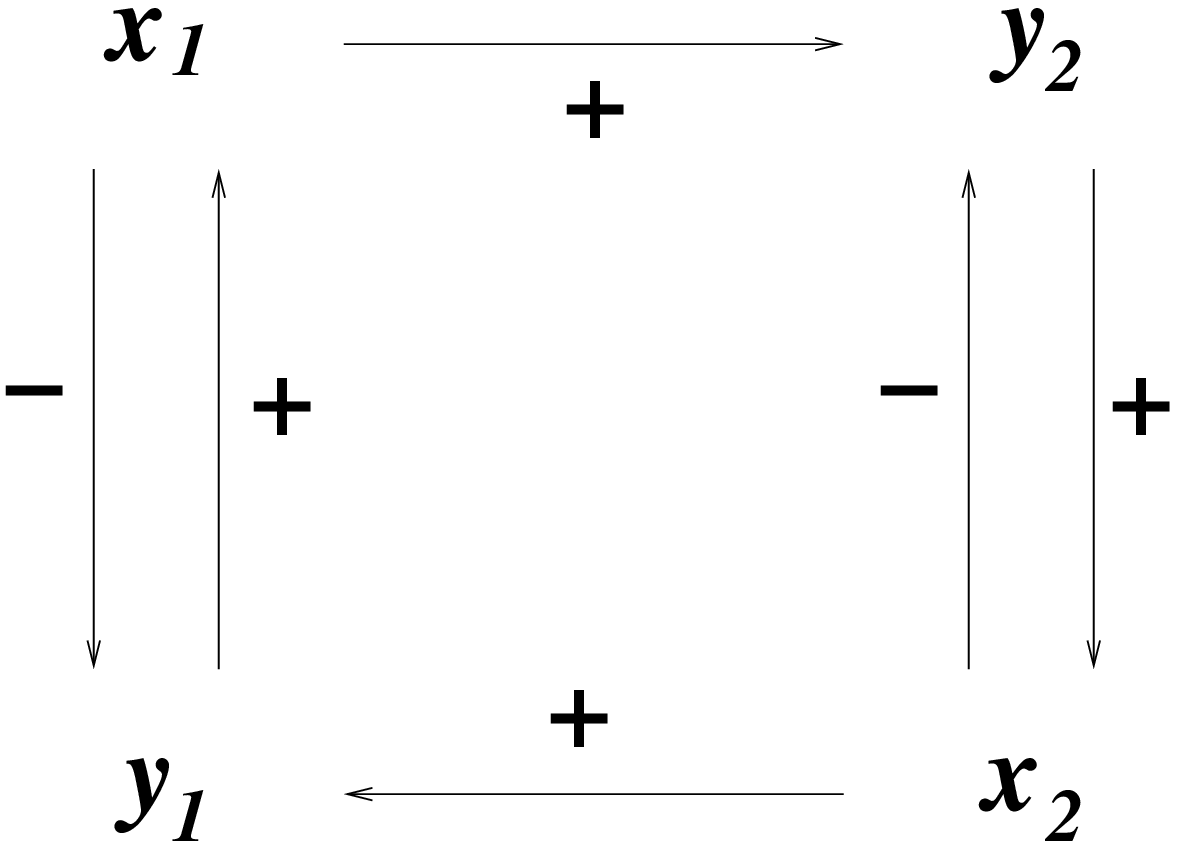}
\caption{Example}
\label{fig:2by2example}
\end{figure}

\medskip
\noindent{\bf Remark:}
The conditions~(\ref{reducedpartialexample}), (\ref{betacond1}),
and~(\ref{betacond2}) are satisfied, in particular, if one assumes the
following easier to check conditions on the functions $\beta_i$'s,
$\alpha_l$'s, and $\gamma_i$'s, 
The functions $\beta_i$ are asked to be so that:
\[
\frac{\partial \beta_i}{\partial x_k}(x)<0
\]
for every $i,k=1, \dots, n$ such that $i \not=k$ (cooperativity condition
among $x_i$ variables), and also so that:
\begin{eqnarray}
\lim_{x_1 \rightarrow +\infty, \dots, x_n \rightarrow +\infty}
 \beta_i(x_1,\dots,x_n)=+\infty 
 \label{cond:infty1}
 \\
\lim_{x_1 \rightarrow -\infty, \dots, x_n \rightarrow -\infty} \beta_i(x_1,\dots,x_n)=-\infty.
 \label{cond:infty2}
\end{eqnarray}
These last conditions are very natural.  They are satisfied, for example, if
there is a linear decay term $-\lambda_ix_i$ in the differential equation for
each $x_i$, and all other variables appear saturated in this rate.
Since $\frac{\partial \beta_i}{\partial x_k}(x)<0$ for all $i \not=k$,
(\ref{cond:infty1})-(\ref{cond:infty2}) imply that conditions~(\ref{betacond1})
and~(\ref{betacond2}) both hold.
Regarding the remaining functions, we ask:
\[
\sum_{l=1}^m 
\frac{\partial \gamma_i}{\partial y_l}\frac{\partial \alpha_l}{\partial x_k}
\leq 0
\]
for all $i,k=1, \dots, n$ such that $i \not=k$.
This condition can be guaranteed to hold based only upon the signs of
the partial derivatives: it holds true if there is no indirect negative
effect (through the variables $y_l$) of any variable $x_k$ on any other
variable $x_i$.
The diagram shown in Figure~\ref{fig:2by2example} illustrates one such
influence graph (signs indicate signs of partial derivatives), for $n=m=2$.
Observe that this example cannot describe a monotone system (with respect to
any orthant cone, i.e., it is not cooperative under any possible change of
coordinates of the type $x_i\rightarrow-x_i$ or $y_l\rightarrow-y_l$).
An entirely analogous example can be done for any $n=m$, the key
property being that each variable $x_i$ ``represses'' its associated variables
$y_i$ and the $y_l$'s ``enhance'' some or all other variables.

\section{Appendix}

\begin{defn}
The closed half-space $H^l \subset \mathbb R^l$, is defined as follows:
\[
H^l=\{(x_1,x_2, \cdots, x_l) \in \mathbb R^l \ | \ x_1 \geq 0 \}.
\]
The boundary of $H^l$, denoted as $\partial H^l$, is $\mathbb R^{l-1}$.
\label{def:half_space}
\end{defn}
\begin{defn}
A subset $M \subseteq \mathbb R^N$ is called a $l$-dimensional $C^r$ manifold with boundary if it satisfies:
\begin{enumerate}
\item There exists a countable collection of open sets $V^\alpha \subseteq
  \mathbb R^N$, $\alpha \in \mathcal I$, where $\mathcal I$ is some countable
  index set, so that, with $U^\alpha \equiv V^\alpha \bigcap M$, one has
  $M=\bigcup_{\alpha \in \mathcal I} U^\alpha$. 
\item There exists a $C^r$ diffeomorphism $x^\alpha$ defined on each
  $U^\alpha$ which maps $U^\alpha$ to some set $W \bigcap H^l$ where $W$ is
  some open set in $\mathbb R^l$. 
\end{enumerate}
The boundary of $M$, denoted as $\partial M$, is the set of points in $M$ that
are mapped to $\partial H^l$ under $x^\alpha$, for some $\alpha \in \mathcal
I$.  
%formatted-06posta.tex
\label{def:manifd_bdy}
\end{defn}
A compact manifold is a manifold that is compact as a topological space. The definition implies that $\partial M$ is a well defined $C^r$ manifold of dimension $l-1$, and $M \setminus \partial M$ is an $l$-dimensional $C^r$ manifold; see \cite{Milnor} for the details.

\end{document}